\documentclass{article}
\usepackage{arxiv}
\usepackage[utf8]{inputenc}
\usepackage[T1]{fontenc}    
\usepackage{hyperref}      
\usepackage{url}           
\usepackage{booktabs}      
\usepackage{amsfonts}
\usepackage{nicefrac}
\usepackage{microtype}
\usepackage{lipsum}
\usepackage{graphicx}
\graphicspath{ {./images/} }
\usepackage{comment}
\usepackage{xcolor}
\usepackage{amssymb}
\usepackage{stix}
\usepackage{enumitem}
\usepackage{amsmath}
\usepackage{relsize}
\usepackage{cite}
\title{Axioms for the Measure of Evidence}

\author{
 Christopher D. Fiorillo \\
  Department of Bio and Brain Engineering, KAIST\\
  Daejeon, South Korea\\
  \texttt{rationalobserver137 at gmail} \\
\and
 Min Sheo Choi \\
  Holistic Manifold Inc.\\
  Daejeon, South Korea\\
  \texttt{choi at holisticmanifold dot com} \\
\and
 Jaime Gómez-Ramírez \\
  Institute of Biomedical Research and Innovation\\
  University of Cádiz, Puerta del Mar Hospital\\
  Cádiz, Spain\\
  \texttt{borriquear at gmail} \\
 }

\begin{document}
\maketitle
\begin{abstract}
There has not been an established mathematical measure of evidence.  Some Bayesians have argued that probability can be an objectively correct measure of ``\textit{rational degrees of belief},'' which we do not distinguish from evidence.  However, support for the objectivist view has been limited due to the lack of a general method for assigning probabilities to evidence (belief) \textit{de novo}.  The standard axioms of Kolmogorov and Cox solve only \textit{the calculation problem}, specifying how probabilities can be calculated from other probabilities. They do not solve \textit{the measurement problem} of how to determine the uniquely correct value of $P(A)$ given only $A$.  The prototypical solution has always been Laplace's \textit{principle of indifference}, which assigns equal (uniform) probabilities to all possibilities.  However, uniformity is well known to be inconsistent with the standard axioms when there are infinite possibilities. Here we introduce new axioms that resolve this inconsistency.  We first show that all measures must ultimately be based on uniformity, and that uniformity is inevitable if all propositions are adequately defined.  We then reconcile uniformity with infinite possibilities by using hyperrational numbers, so that an infinite sum of infinitesimal probabilities can equal one.  Our axioms thereby provide a general and relatively simple solution to the measurement problem.  We discuss a variety of conceptual obstacles that have made this solution difficult to recognize. 
\end{abstract}


\section{\label{sec:introduction}Introduction}
Classical logic has proven remarkably useful, despite classifying propositions in a binary manner as either certain to be false or certain to be true (0 or 1). The propositions of most interest are those that may or may not be true given our incomplete knowledge.  We have some amount of evidence that a proposition is true, and our evidence can vary between the extremes of certain falsehood and certain truth.  Therefore we need a graded measure of evidence. 

We propose that probability is such a measure. There has been no consensus about what probability measures, but it is not commonly thought to be a mathematical measure of evidence \cite{Berger1985,bernardo1994,Jaynes2003,beisbart2011}. Bayesians view probability as a measure of \textit{belief} (or credence or plausibility). \textit{Belief} has the connotation of being aphysical, weakly justified, and possibly irrational. We view this as unfortunate, since we see human beliefs as corresponding to internal brain states that are evidence, however weak, for external and future states \cite{Fiorillo2008,Fiorillo2012}. Regardless of their physical basis, we see no useful mathematical distinction between \textit{evidence} and \textit{belief}.  We use them as synonyms but strongly prefer \textit{evidence}. 

There has not been consensus among Bayesians as to whether probabilities can be objectively correct measures of evidence (belief).  All Bayesians recognize that evidence and probability are subjective in the sense that they vary from one observer to another (they are local in space and time). Only objectivists believe that there is a single uniquely correct number $P(A|B)$ that is fully determined by reason alone given only the definition of propositions $A$ and $B$ \cite{Jeffreys1961,Cox1961,Jaynes2003,Laplace1814}.  We can formally express this ideal for $P(A)$ by stating that there exists at least one set of axioms $\mathcal{X}$ such that  
\begin{equation}
\exists \mathcal{X},\; A {\implies} P(A)\; \forall A
\label{eq:A_implies_P(A)}
\end{equation}

Determining $P(A)$ from $A$ is what we call \textit{the measurement problem}. The standard axioms of Kolmogorov (K-axioms) and Cox (C-axioms) address only \textit{the calculation problem}, specifying how to calculate one probability from another \cite{kolmogorov1950,Cox1961}.  Axioms $\mathcal{X}$ would provide a general and universal solution for both measurement and calculation, but even their existence has been controversial. We propose the present axioms as an instance of $\mathcal{X}$.  

Though everyone sees objectivity as desirable, we believe most Bayesians to be ``subjectivists.'' They see the relation of $P(A)$ to $A$ as at least partially indeterminate in many if not all cases, so that two observers need not agree about probabilities, even given the same information and propositions \cite{Berger1985,bernardo1994,Kass1996,Lindley2000,Howson2003}.  Indeed, the relation \textit{is} indeterminate and ambiguous, at least as a practical matter, so long as we lack a solution to the measurement problem.  A general solution may help to overcome the longstanding skepticism of subjectivists.  

Objectivists have proposed various solutions. These include Laplace's \textit{principle of indifference}, Jeffreys's \textit{invariance}, and Jaynes's \textit{maximum entropy}.  However, all of these have been controversial in one respect or another, and none are seen as entirely general  \cite{Berger1985,bernardo1994,Kass1996,Jaynes2003,Robert2009}.  The principle of indifference (uniformity) simply says that in the absence of evidence favoring one proposition over another, the probability of each must be equal.  This is pure reason, and it must be true if probability measures evidence.  The justification for invariance and maximum entropy is less obvious, and they have been only loosely grounded in propositional logic.  

The specific and well known problem is that none of these methods are consistent with either the K-axioms or C-axioms \cite{kolmogorov1950,Cox1961} when applied to infinite sets of propositions in the absence of prior empirical evidence \cite{Berger1985,bernardo1994,Kass1996,Jaynes2003}. This may appear to be a special case of limited significance. However, the set of all possibilities is infinite if we are to be realistic.  Furthermore, a universal and objectively correct measure of evidence should not assume the truth of any empirical observation (\ref{sec:non-empirical}).

We draw attention to four distinctive features of our approach.  First, we build on the work of Jeffreys by insisting that empirical knowledge is irrelevant to $P(A)$ except insofar as it is inherent to the definition of $A$. Second, we demonstrate that uniformity must underlie all probabilities. Skeptics have emphasized that uniformity over one set of propositions (and parameters) is typically inconsistent with uniformity over another, but we argue that challenge can be solved (\ref{sec:disc_multiple_parameters}). Third, we overcome the conflict between uniformity and infinite sets by introducing a new measure of evidence $E(A)$ that can be infinite or infinitesimal.  Fourth, we define probability as the fraction $P(A) \equiv E(A) / E(A{\vee}\neg A)$. 

Our axioms support the original and intuitive ``classical definition'' of probability advocated by Laplace over $200$ years ago \cite{Laplace1814}.  The substantial challenge is to explain why the classical framework is sufficient as a technical matter despite apparent shortcomings that have been debated for over a century.  Therefore much of this manuscript is devoted to overcoming conceptual obstacles.     

\section{\label{sec:desiderata}Desiderata}
Here we list ten desirable features of a mathematical measure of evidence, followed by three specific mathematical constraints our axioms should satisfy.  We believe our axioms satisfy all of these criteria.

\subsection{\label{sec:desiderata_general}General Desiderata}
\begin{enumerate}
    \item \textit{Objective and Universal}: Evidence is subjective insofar as it is local in space and time, but we should have measures of evidence that are objectively correct and universal (the same for all observers).
    \item \textit{General}:  We should be able to measure the evidence for any adequately defined propositions $A,B,\ldots$.
    \item \textit{An Elementary Measure}: We should have a primitive measure of evidence $E(A)$ that is analogous to a distance metric in geometry.
    \item \textit{A Relative Measure}: We should have a measure such as $P(A)$ that compares the evidence for $A$ to its negation.
    \item \textit{A Relational Measure}: We should have a measure such as $P(A|B)$ of the evidence in $B$ for $A$. 
    \item \textit{Quantifiable}: We should have general rules for assigning a unique number to a relative measure such as $P(A)$. 
    \item \textit{Realistic}: Reality presents us with an infinite set of possibilities, and our measures should reflect that.  
    \item \textit{A Calculus}:  We should have rules for calculating one measure of evidence from another. 
    \item \textit{Exact without Assumptions}:  Our measures should be exact in the absence of questionable assumptions.
    \item \textit{Intuitive}: We would hope that our measures can be intuitive, like uniformity.
\end{enumerate} 

Of these, only \textit{the calculation problem} (Desideratum $8$) has a widely recognized and accepted solution, provided by the K-axioms and C-axioms.  

\subsection{\label{sec:desiderata_specific}Specific Desiderata}
Here we list three specific desiderata that serve as mathematical constraints that our axioms should satisfy.  We highlight them because, taken together, they conflict with the standard probability axioms and their application to measure theory.  

\begin{enumerate}
    \item \textit{Infinity}: The set of possible realities (that are conceivable and logically consistent) is infinite.
    \item \textit{Uniformity}:  Uniformity underlies all mathematical measures, including the measure of evidence.
    \item \textit{Impossibility}:  The evidence for a possibility must be greater than the evidence for an impossibility.
\end{enumerate} 

The K-axioms are not consistent with the combination of $1$ and $2$. They assert that probability is a real number, and that the sum of all probabilities must equal one \cite{kolmogorov1950}.  Yet there is no real number that, when summed infinitely, equals one ($\forall r{\in}\mathbb{R}, r+r+\ldots \neq 1$).

Measure theory has generally complied with Specific Desiderata $1$ and $2$, but has contradicted $3$ by assigning measure zero to all ``null sets,'' which includes both the logically impossible (the empty set \eqref{eq:false}) and all finite subsets of possibilities \cite{Kapinski2004}.  That works as an excellent approximation for many practical purposes, but it cannot be correct to assign probability zero to both the possible and impossible.   

\section{\label{sec:non-empirical}Segregation of Empirical and Non-Empirical Evidence}
We propose below that $P(A|B)$ measures the evidence in any proposition $B$ for any proposition $A$.  We naturally have a particular interest in cases in which the truth of $B$ is known from observation, and therefore $P(A|B)$ can be said to measure empirical evidence $B$ for $A$.  However, this is not the most general way to understand $P(A|B)$, since it is the same number regardless of any knowledge we may have about the truth of $A$ or $B$ or $AB$.\footnote{We regard `knowledge' and `information' and `evidence' to be synonymous conceptually, and physically, but distinguished mathematically by the fact that evidence \textit{for} $A$ has a polarity that information/knowledge \textit{about} $A$ does not.} It is better understood as a measure of non-empirical evidence, meaning evidence from reason alone.  

It may seem that all evidence is empirical, but this is not the case.  For example, reason alone indicates that, for any propositions $A$ and $B$, the evidence that both are true must be less than or equal to the evidence that $A$ is true, $E(A {\land} B) \leq E(A)$.  Equality applies only in the rare and rather trivial case that we know that $B$ must be true if $A$ is true. Otherwise $E(A {\land} B) < E(A)$.  This is obvious without even specifying what $A$ and $B$ represent, rendering the problem entirely free of empirical evidence. This demonstrates that evidence can be non-empirical, a consequence of pure reason.  

Our axioms enforce segregation such that all that is empirical in nature must be expressed in propositions, and probabilities express only rational (non-empirical) relations between propositions.  This has two major advantages.  First, this segregation is essential if our measure is to be objective and universal (\ref{sec:desiderata}). Empirical evidence is local in space and time, and thus subjective in the sense that it varies across observers. However, because our probabilities measure only non-empirical evidence, two observers using our axioms properly will never disagree about a probability.  One might have empirical knowledge $A$ and another $B$, so that they disagree about the truth of $A$ and $B$, and yet both will agree about the numbers assigned to $P(A)$, $P(A|B)$, etc.

A second major advantage of measuring only non-empirical evidence is that it greatly simplifies the measurement problem.  A basic challenge is that we tend to have a lot of empirical evidence from diverse sources. It is usually the case that the more we know, the more difficult it is to quantify, because our knowledge breaks symmetries.  Reasoning without empirical knowledge provides greater symmetry and simpler mathematics.   

Jeffreys proposed that we should first assign numbers to probabilities starting from the complete absence of information \cite{Jeffreys1932,Jeffreys1961,Jaynes2003}.  The meaning of ``complete absence'' has been controversial \cite{Irony1997,Jaynes2003}. In our view it means using no empirical knowledge of $A$ other than that within the definition of $A$.  Therefore the relation between $P(A)$ and $A$ is pure reason, and entirely non-empirical.

Reasoning from complete agnosticism is unnatural.  It requires that we use our imagination to purge our minds so that we see nothing as fact and everything as possible (as though we are God before Creation and nothing exists).  This purge allows us, at least in principle, to achieve what Nagel called ``\textit{the view from nowhere}'' \cite{Nagel1986}.  It is not a denial of local and subjective knowledge, but rather the means to objectively quantify it.   

\section{\label{sec:logic}Logical Identities}
$A,B,C,\ldots$ symbolize generic propositions.  We specify all logical relations using negation ($\neg A \equiv not\ A$), disjunction ($A {\vee} \neg A, \vee \equiv OR$), and conjunction ($AB \equiv A{\land}B,\ \land \equiv AND$).  These give rise to various logical identities, such as $\neg(AB) = (A\neg B) \vee (\neg AB) \vee (\neg A \neg B)$.  Below we highlight three that may not be obvious.

The first expresses the fact that a proposition is defined by default to be inclusive rather than exclusive, concerning some aspect of reality while being agnostic about others. In proposing that $A$ is true, we are implying agnosticism about the truth of $B,C, \ldots$, rather than implying that they are false.  $A$ is the same proposition whether this agnosticism is implicit or explicit. Therefore the following three expressions are equivalent.
\begin{equation}
A\; = A \land (B \vee \neg B)\; = AB \vee A \neg B
\label{eq:inclusive}
\end{equation}

$AB \vee A\neg B$ conveys more to us than $A$, by expressing $A$ in the larger context of $B$, but it is nonetheless the same proposition.  Obviously we could extend this to explicitly express our agnosticism about any number of additional propositions $C,D,E,\ldots$. It can be advantageous to choose a number, but the choice is substantially arbitrary and will not matter with respect to $P(A)$ (\ref{sec:model}, \ref{sec:th_scale}).  
From the above identity it follows that
\begin{equation} 
\top\; \equiv\; A \vee \neg A\; =\; B \vee \neg B
\label{eq:true}
\end{equation}

\noindent 
where `$\top$' is a proposition that must be true as a matter of logic alone. That $A {\vee} \neg A\; =\; B {\vee} \neg B$ can be seen by replacing $A$ with $A{\vee}\neg A$ in the previous identity \eqref{eq:inclusive}, and considering an analogous equation with $A$ and $B$ exchanged.  Then one can show that both $A{\vee}\neg A$ and $B{\vee}\neg B$ are equal to $AB {\vee} \neg AB {\vee} A\neg B {\vee} \neg A\neg B$, which can be expressed as \textit{anything is possible}.  

Its negation, \textit{nothing is possible}, is $\neg(A {\vee} \neg A) = \neg(B{\vee}\neg B)$, meaning that $A$ and $\neg A$ are either both true, or both false.  But these are identical, since $\neg A {\land} \neg \neg A = A {\land} \neg A$.  Therefore $\neg(A {\vee} \neg A) = A {\land} \neg A$, and likewise for $B$.  By negating equation \ref{eq:true} we have
\begin{equation}
\bot\; \equiv\; A \land \neg A\; =\; B \land \neg B
\label{eq:false}
\end{equation}

\noindent
where $\bot$ is defined as a proposition that is false because it is a logical contradiction and impossibility. It also corresponds to the empty set, $\bot {\equiv} \{\}$, since it is the negation of all possibilities in an exhaustive set.

\section{\label{sec:possibility}The Space of Possibilities}
Before assigning probabilities we must determine what is possible by precisely defining all propositions. We argue that this is the only major challenge, but it tends to be substantially more difficult than it appears (\ref{sec:discussion}). 

We are free to define our propositions however we like.  But if we want consistency between our beliefs, our definitions, and our measures, all of our relevant empirical knowledge must be formally and precisely defined in our propositions (\ref{sec:non-empirical}).  For example, if we only distinguish two possibilities, $A {\vee} \neg A$, then $P(A) {=} 1/2$, regardless of what $A$ is (e.g., \textit{I exist}, \textit{God exists}, \textit{the sun will rise tomorrow}, etc.).  That may satisfy no one, but there is no problem here with the mathematics. Fault can only lie in the fact that we chose to define $A$ and $\neg A$ as each being just one indivisible possibility (\ref{sec:uniformity}).

\subsection{\label{sec:model}$\mathbb{U}$ is a model of what is possible}
We can never be certain of what is possible, since our concern is usually the future and possibilities are inherently unobservable. But we can deduce a set of possibilities from a belief about what exists in reality.  Therefore the first step should be to formalize beliefs by specifying a precise ontological model \textbf{O}.  For our statistical model to be as general as possible, \textbf{O} would ideally be a detailed model of every universe we can imagine.  But \textbf{O} could be as simple as `\textit{there exists a coin},' or `\textit{there exists a standard deck of $52$ cards}' (though these terms must be unambiguously defined). By assuming \textbf{O} is true we can logically deduce the exhaustive set $\mathbb{U}$ of all mutually exclusive possibilities
\begin{equation}
\mathbf{O} \implies \mathbb{U} = \{u_1,u_2,\ldots\}
\label{eq:O}
\end{equation}

\noindent
where $\mathbb{U}$ could be a finite or infinite set depending on \textbf{O}, and the subscripts are arbitrary labels that are not intended to imply an order. That $\mathbb{U}$ should be deduced from a model may seem obvious, but we highlight it because deduction can be challenging (\ref{sec:disc_multiple_parameters}), and because there is typically no attempt at deduction in problems with continuous variables (\ref{sec:disc_absolute_space}, \ref{sec:disc_numeric}). 

Each possibility $u_i{\in}\mathbb{U}$ is indivisible, and is therefore known as an ``atomic proposition.''  If anyone disputes that $u_i$ is indivisible, this merely indicates that they prefer an alternative model \textbf{O}.  We can consider multiple models, but a particular model \textbf{O} must be unambiguously defined at the outset, and $u_i$ must be indivisible ($\forall i$). This is necessary for propositions to be well defined, and it insures uniformity (\ref{sec:uniformity}).  

However we choose to define \textbf{O}, $\mathbb{U}$ becomes our entire universe of possibilities.  Model \textbf{O} is itself a proposition (e.g., perhaps the coin does not exist, or perhaps it has more than two sides).  But since $\mathbb{U}$ assumes that \textbf{O} must be true, $\neg \mathbb{U}$ and $\neg \mathbf{O}$ correspond to logical contradiction (the empty set) \eqref{eq:false} and zero evidence. A useful measure of the evidence for \textbf{O} would require a new model that includes a plausible alternative to \textbf{O}. \textbf{O} would ideally be a model so comprehensive that it includes every possibility (and thus every model) we can conceive.  

\subsection{\label{sec:spaces}Three Spaces}
We refer to any exhaustive set of mutually exclusive possibilities as a `\textit{space}.'  To aid understanding it is helpful to distinguish three types that we denote $\top$-space, U-space, and S-space, each with a corresponding set $\mathbb{T},\mathbb{U},\mathbb{S}$.  $\top$-space is the true and real space of all that ever has or ever will be possible in our universe.  It is the only space that is truly exhaustive and unfalsifiable, since it is the only space that is not associated with a model. However, we cannot know what it is, and therefore it is of no direct relevance to mathematics. We mention $\top$-space only because it is critical that we do not mistake $\mathbb{U}$ for $\mathbb{T}$.

The `U' in ``U-space'' signifies uniform, unit, and universe. $\mathbb{U}$ is the entire universe of possibilities within our model \textbf{O}, each possibility $u_i{\in}\mathbb{U}$ is a fundamental and indivisible unit (\ref{sec:uniformity}), and evidence is uniformly distributed over these units (\ref{sec:uniformity}). 

$\mathbb{U}$ in our axioms corresponds to the ``sample space'' $\Omega$ in the K-axioms \cite{kolmogorov1950}. All propositions (``events'' in the terminology of the K-axioms) are subsets of $\mathbb{U}$, with the powerset $\mathcal{P}(\mathbb{U})$ being the set of all propositions.

However, $\mathbb{U}$ and $\Omega$ usually differ in practice, since evidence can have any distribution over $\Omega$ but will be uniform over $\mathbb{U}$ in all cases (\ref{sec:uniformity}). Non-uniformity arises only when we have propositions composed of differing numbers of elementary units $u_i$. For example, $\mathbb{U}$ could be all combinations of the $52$ cards in a standard deck.  From these we can create many additional ``state spaces'' (S-spaces), $\mathbb{S}_1,\mathbb{S}_2,\ldots$, each of which is a distinct grouping of the $52$ cards into disjoint and exhaustive subsets.  For example, a specific S-space would consist of $3$ subsets (states) corresponding to aces, face cards, and cards numbered $2$ to $10$.  By counting the number of possibilities we find the probability that an unknown card is in one of these groups to be $4/52$, $12/52$, and $36/52$, respectively. 

In general, there will be many S-spaces defined on a single U-space. The evidence will be non-uniformly distributed for all except the one special case $\mathbb{S}_i {=} \mathbb{U}$.  We distinguish S-spaces to show how non-uniformity arises, and to emphasize that what is commonly known in practice as a ``state space''  or ``sample space'' $\Omega$ usually does not correspond to a U-space $\mathbb{U}$.  

\section{\label{sec:uniformity}Uniformity}
Uniformity has always been the primary solution to the measurement problem, but it has been widely criticized \cite{Berger1985,bernardo1994,Kass1996}. Even the most committed of objective Bayesians have not viewed it as a universal solution \cite{Jeffreys1961,Jaynes2003}.  We claim that it is as universal and incontrovertible as mathematics itself. Here we present two rationales for why uniformity must underlie all probabilities.     

There are various definitions of a natural number, all based on sets (\ref{sec:disc_numeric}).  Uniformity is essential to all, since each element of a set is counted equally, regardless of the distinctions between the elements.  We define a number as the measure (size or cardinality) of a set.  Thus $3$ is the measure $\mu$ of any set that contains three elements, 
\begin{equation}
3\; \equiv\; \mu\{\bullet,\bullet,\bullet\} =\; \mu\{X,{\scriptstyle X},{\scriptscriptstyle X}\}\; =\; \mu\{*,\#,@\}\; =\; \ldots
\label{eq:three}
\end{equation}

\noindent
such as a set composed of two mice and one elephant.  They are each counted equally because they equally share some common feature that justifies their inclusion as elements of one set, like all being mammals.  

Since the natural numbers are based on uniformity, and all numbers are derived from the natural numbers, all numbers are based on uniformity.  Therefore any numbers that measure evidence must also be based on uniformity.

An apparently distinct argument for uniformity is more specific to evidence.  Since we are measuring only the non-empirical evidence for proposition $A{\subseteq}\mathbb{U}$ (\ref{sec:non-empirical}, \ref{sec:possibility}), the measure must depend only on the definition of $A$. $A$ is adequately defined only if \textit{every possibility} that would render it true or false is specified (\ref{sec:spaces}).  Therefore a proposition must be a disjunction of indivisible ``atomic propositions,'' as in $A = u_1{\vee}u_2$ and $\neg A = u_3{\vee}u_4{\vee}u_5$.  Reason alone can never justify the assignment of unequal evidence to indivisible and mutually exclusive propositions.

Uniformity has always been implicit within all use of numbers for measurement. Its validity cannot be proven because it is inherent to reason itself, and therefore a necessary premise of proofs involving measures. The utility of a proof is in proving a non-obvious conclusion from an obvious premise.  Uniformity is too obvious to prove. Therefore we assert uniformity in our first axiom (albeit implicitly).

\section{\label{sec:hypernatural}Infinite Sets}
The set of real possibilities, given our limited knowledge, is infinite.  This is seen with the continuum, but it is also implicit in the agnostic and inclusive nature of propositions \eqref{eq:inclusive}.  Therefore we need the capacity to measure infinite sets of differing sizes.  The inability to do this with real numbers underlies the incompatibility of the K-axioms and C-axioms with uniformity over infinite sets.   

A related problem arose historically with the infinitesimals of calculus.  The standard solution has been to use limits, although this causes confusion about the status of `$dx$.'  An alternative known as ``Non-standard Analysis'' solved the problem by introducing ``hyperreal numbers'' \cite{Robinson1966}.  These extend the real numbers to include infinite and infinitesimal numbers that follow the standard rules of arithmetic.  
For present purposes we need only hypernatural and hyperrational numbers \cite{Lovyagin2021}.  This provides a simple resolution to the measurement problem.  However, it creates tension with conventional views of the continuum and the cardinality of infinite sets, as we discuss further below (\ref{sec:disc_countably_infinite}). 

\section{\label{sec:axioms}Axioms}
\subsection{\label{sec:ax_measure}The Measure of a Set}
The measurement problem of probability theory has been debated as though it is specific to probabilities, but we see it as concerning the foundation of mathematics, and particularly the definition of numbers.  We therefore begin by defining every hypernatural number to be the measure $\mu$ (cardinality) of any set (finite or infinite) containing the corresponding number of elements.  Thus, $\forall \mathfrak{n}{\in}\mathbb{N}_H$,
\begin{flalign}
\textsf{Axiom 1} && \mathfrak{n} \equiv \mu\{x_1,x_2,\ldots x_{\mathfrak{n}}\} &&
\label{eq:ax_uniform}
\end{flalign}

\noindent
where the indices are arbitrary and exchangeable labels that should not be interpreted to imply an order. This defines $\mathfrak{n}$ as the familiar ``counting measure'' based on uniformity, so that $\forall i,j,\; \mu(x_i) = \mu(x_j) = 1$.  However, as natural as this is, it is only a convenient assignment, since there can be no uniquely correct choice of a number for an ``absolute scale'' (\ref{sec:ax_scale}) or for the cardinality of a set of possibilities (\ref{sec:logic}).  Our ultimate concern is only ratios of measures, and any chosen scale will cancel (\ref{sec:th_scale}).  

\subsection{\label{sec:ax_prop_is_set}A Proposition is a Set of Possibilities}
An ontological model logically implies a set $\mathbb{U}$ of mutually exclusive possibilities \eqref{eq:O}. A proposition is a subset $A{\subseteq} \mathbb{U}$.
\begin{flalign}
\textsf{Axiom 2} && A \equiv \{u_1 \vee u_2 \vee \ldots u_{\mathfrak{n}}\}&&
\label{eq:ax_prop_is_set}
\end{flalign}

\noindent
As emphasized above (\ref{sec:model}), $A$ is adequately defined only if \textit{every} possibility $u_i$ within $\mathbb{U} = A \cup \neg A$ is specified.  The number of possibilities in $A$ and $\neg A$ must each be specified, but not every possibility must be enumerated or defined in detail (this is unattainable for an infinite set) (\ref{sec:disc_countably_infinite}). 

\subsection{\label{sec:ax_E_is_a_measure}Evidence is a Measure}
Using Axioms 1 and 2, we define the measure of evidence $E(A)$ for $A$ to be the measure of set $A$,   
\begin{flalign}
\textsf{Axiom 3} && E(A) \equiv \mu(A)&&
\label{eq:ax_E_is_a_measure}
\end{flalign}

\noindent
and thus a hypernatural number, $E(A){\in}\mathbb{N}_H$.

\subsection{\label{sec:ax_scale}Optional Axiom for Infinite Scale}
Axiom 4 provides several advantages but is not essential.  It assigns by convention an infinite number of elements $\aleph{\in}\mathbb{N}_H$ ($\aleph{>}n\ \forall n{\in}\mathbb{N}$) to $\mathbb{U}$, so that
\vspace{3pt}
\begin{flalign}
\textsf{Axiom 4}\ && \mathbb{U} \equiv \{u_1,u_2,\ldots u_\aleph\} &&
\label{eq:ax_scale}
\end{flalign}

\noindent
and $|\mathbb{U}| {=} \aleph$.  The particular assignment can only be a matter of convenience.\footnote{It is convenient for $\aleph$ to be evenly divisible by $2^n, \forall n{\in}\mathbb{N}$.  For example, it could be $\aleph {=} 2^{\aleph_0}$, where $\aleph_0 {=} |\mathbb{N}|$ is the number of natural numbers.  Note that even if one chose to equate $\aleph$ with $\aleph_0$, $\aleph$ is not a transfinite number as defined by Cantor.}

Although the utility of Axiom 4 is most obvious when $\mathbb{U}$ must necessarily be an infinite set, it applies to all cases, even those that we would normally regard as finite.  This has several advantages discussed below (\ref{sec:th_scale}).  One is that it more accurately reflects reality.  For example, when we consider \textit{heads} and \textit{tails} as the outcomes for a coin toss, we need not pretend that there are only $2$, or any finite number, of possible states of reality. It is not necessary to specify the infinite number of possibilities, because to find $P(heads)$ we are merely partitioning $\mathbb{U}$ into two subsets, those featuring heads versus tails, and then comparing their relative size.  

\section{\label{sec:theorems}Theorems}
From our axioms we can derive the K-axioms and C-axioms, and thus all the familiar theorems of probability theory.  Our Theorems 1 and 2 can be proven from Axioms 1-3, and Theorem 6 from Axiom 4 (with proofs so simple and concise that we do not highlight them). We can also use $E$, which measures the evidence for a single proposition in isolation, to invent any relative measures we like. Theorems 3-5 define odds, probability, and conditional probability.

\subsection{\label{sec:th_additivity}Additivity}
Since Axiom 1 asserts that the measure of every set is the number of elements in it, measures must be additive such that for any collection of $\mathfrak{n}$ sets $X_1,X_2,\ldots X_{\mathfrak{n}}$ that are disjoint ($X_i \cap X_j = \{\}, \forall i,j$),  
\begin{flalign}
\textsf{Theorem 1} && \mu(\bigcup_{i{=}1}^{\mathfrak{n}} X_i) = \sum_{i{=}1}^{\mathfrak{n}} \mu(X_i) &&
\label{eq:th_additivity}
\end{flalign}

\noindent
This is the standard additivity of measure theory, and a more general form of the third K-axiom \cite{kolmogorov1950}.  However, it is obtained here without the complexity that is introduced by considering sets of real numbers.  We believe that, just as propositions are necessarily discrete and countable, all well defined sets with ontological status should be countable (\ref{sec:disc_countably_infinite},\ref{sec:disc_numeric}).   

\subsection{\label{sec:th_sum_rule}The Sum Rule}
Since $A{\cap}\neg A = \{\}$, and since $\mathbb{U} = A{\vee}\neg A = \top$ is a single proposition \eqref{eq:true}, from equations \ref{eq:ax_prop_is_set}, \ref{eq:ax_E_is_a_measure}, and \ref{eq:th_additivity} we have ``the sum rule.''
\vspace{3pt}
\begin{flalign}
\textsf{Theorem 2} && E(\top) = E(A) + E(\neg A) &&
\label{eq:th_sum_rule}
\end{flalign}
\noindent

This can be used to derive the corresponding rules of the C-axioms and K-axioms (which could be added here as additional theorems) \cite{kolmogorov1950,Cox1961}. It is more fundamental in that it applies to elementary evidence $E$ rather than relative evidence $P$.  Note that according to our definition of $P$ \eqref{eq:th_defineP}, $P(A) + P(\neg A) = 1$ (Cox's sum rule) simply because it is the sum of complimentary fractions.  This eliminates the need for the lengthy justification of Cox's sum rule given previously \cite{Jaynes2003}.  It also eliminates the need for the second K-axiom, which assigns $1$ by convention \cite{kolmogorov1950}.

\subsection{\label{sec:th_odds}Odds is a Relative Measure}
We define odds $O(A)$ as the ratio of the evidence that $A$ is true to the evidence that it is false.
\begin{flalign}
\textsf{Theorem 3} && O(A) \equiv \frac{E(A)}{E(\neg A)} &&
\label{eq:th_odds}
\end{flalign}

\noindent
Since it follows that $O(\neg A) \equiv E(\neg A)/E(A)$, we also have $O(A)\ O(\neg A) = 1$, which is analogous to Cox's sum rule for probability.  We can convert the product to summation by defining log-odds $L(A) \equiv log(O(A)) \equiv log(E(A)) - log(E(\neg A))$, so that $L(A) + L(\neg A) = 0$.  Odds and log-odds exhibit convenient symmetries that are lacking in a fraction like probability. 

\subsection{\label{sec:th_defineP}Probability is a Relative Measure}
We define probability as the following fraction.
\begin{flalign}
\textsf{Theorem 4} && P(A) \equiv \frac{E(A)}{E(A \vee \neg A)} &&
\label{eq:th_defineP}
\end{flalign}

\noindent
$P(A)$ is a relative measure because it compares the evidence for two propositions.  Although potentially misleading, it can be considered ``unconditional'' insofar as it concerns only a single freely chosen proposition (since the definition of $A$ determines that of $\neg A$).

\subsection{\label{sec:th_definePcond}Conditional Probability is a Relational Measure}
We define conditional probability to be
\begin{flalign}
\textsf{Theorem 5} && P(A|B) \equiv \frac{E(AB)}{E(B)}&&
\label{eq:th_definePcond}
\end{flalign}

\noindent
This measures evidence $B$ for $A$. $B$ can be described as the \textit{reference proposition}, and $A$ the \textit{proposition of interest}.  An unconditional probability $P(A)$ is the special case of $P(A{|}B)$ where $B = A{\vee} \neg A$.  We refer to $P(A{|}B)$ as a \textit{relational} measure (\ref{sec:desiderata}) to distinguish it from \textit{relative} measures such as $P(A)$ \eqref{eq:th_defineP} and $P(A)/P(B) = E(A)/E(B)$.  The latter do not depend on the relation (interaction) implicit in the conjunction $A{\land}B$. 

From our definitions of probability and conditional probability (\ref{eq:th_defineP}, \ref{eq:th_definePcond}), and logical identities (\ref{sec:logic}), one can derive $P(A{|}B) = P(AB) / P(B)$ and other expressions of Cox's product rule (a.k.a. Bayes's Theorem) \cite{Cox1961,Jaynes2003}. As for Cox's sum rule (\ref{sec:th_sum_rule}), this justification is considerably simpler than that given previously \cite{Jaynes2003}.

\subsection{\label{sec:th_scale}Optional Theorem for Infinite Scale}
Axiom 4 asserts that $\mathbb{U}$ contains an infinite number of elements, $\aleph{\in}\mathbb{N}_H$.  Since proposition $\mathbb{U} = A{\vee}\neg A = B{\vee}\neg B = \ldots$ is an instance of $\top$ \eqref{eq:true}, it follows from Axioms 1-4 that  
\begin{flalign}
\textsf{Theorem 6}\ && E(\top) \equiv \aleph &&
\label{eq:th_scale}
\end{flalign}
\noindent

Axiom 4 and Theorem 6 are obviously not needed if $\mathbb{U}$ can be defined as finite.  Without them we can follow the standard procedure of assigning a number to $E(\top) {=} |\mathbb{U}|$ on a case by case basis, such $6$ for a standard die. They are also not needed for infinite sets if we do not insist on assigning a number to the evidence $E$ for a single proposition in isolation. Indeed, the axioms of geometry provide no guidelines for assigning a number to a single distance. Aside from convenience, a number can only be assigned in an entirely arbitrary manner. Uniquely correct numbers only arise once we consider ratios of distances, and the same applies to evidence.  Axioms 1-3 are sufficient to determine these ratios.  

An advantage of Axiom 4 and Theorem 6 is in providing a single universal convention that preserves the counting measure. For example, $\mathbb{U}$ has size $\aleph$ even in the familiar case of a coin toss (\ref{sec:ax_scale}).  The evidence for \textit{heads} and \textit{tails} will each be $E {=} \aleph/2$. Dividing this by $E(\top) {=} \aleph$ yields $P {=} 1/2$ \eqref{eq:th_defineP}. Thus we obtain the familiar probability without the unrealistic assumption of a universe of only two states.  

Yet another advantage is didactic. Kolmogorov addressed the same issue by assigning $E(\top) {\equiv} 1$ (his second axiom expressed in our notation) \cite{kolmogorov1950}. Were we to adopt his convention, it would make no difference to relative measures such as probability, but it would sow misunderstanding.  Probability is a ratio of the evidence for one proposition relative to another \eqref{eq:th_defineP}, but by making $1$ the denominator, the second K-axiom promotes the common misconception that probability concerns only a single proposition.  An analogous misconception would occur in physics if we always insisted that distances be measured relative to $1$ meter, and also chose to write `$x$' rather than `\textit{x meters}' (or $x/1$).  

\section{\label{sec:discussion}Discussion}
We believe that our axioms address all our desiderata (\ref{sec:desiderata}).  In particular, they overcome a deficiency of the K-axioms and C-axioms by allowing uniformity to be the general solution to the measurement problem. As a consequence, the only challenge in assigning probabilities is in specifying the set of possibilities and then calculating the relative ``volume'' associated with a proposition.  However, there are multiple conceptual obstacles that need to be overcome before our measure of evidence can be adequately understood and utilized.  Below we discuss obstacles associated with infinite sets, multiple parameters, absolute space, and reification of numbers.  We end with examples of how uniform spaces over continuous parameters can be identified.    
\subsection{\label{sec:disc_infinite}Infinite Sets and the Continuum}
It is well known that the K-axioms are inconsistent with uniformity over an infinite set (\ref{sec:desiderata_specific}) \cite{Berger1985,bernardo1994,Kass1996,Jaynes2003}.  The predominate response has been to question uniformity or infinity or both.  

\subsubsection{\label{sec:disc_past_responeses}Past Responses to the Problem}
\vspace{-2mm}
It has been suggested that infinity is strictly a mathematical concept that does not exist in reality \cite{Jaynes2003}. We agree that the infinite does not exist as an actuality. But the issue here is that which could possibly exist, and we believe it is obvious that this set of possibilities is literally infinite.

Another response to the problem has been to point out that the state of ``complete ignorance'' does not exist in reality, and that a real state of knowledge will be quantified by a conditional distribution that is either non-uniform over an infinite set, or uniform over a finite set \cite{Irony1997}.  We agree.  However, a state of complete ignorance does exist in our imagination as a mathematical object, as does zero and the empty set.  We need this state at the foundation for the same reason we need zero and the empty set \cite{Jaynes2003}. 

Yet another response has been to cast doubt on uniformity, and more generally on the feasibility of objective quantification of evidence (belief) \cite{Berger1985,bernardo1994,Kass1996}.  We have pointed out that uniformity is implicit in any use of numbers for measurement (\ref{sec:uniformity}).

\subsubsection{\label{sec:disc_infinite_solution}Our Solution for Infinite Sets}
\vspace{-2mm}
Our axioms allow for uniformity over infinite sets by defining $P(A)$ to be a ratio of more primitive measures that are hypernatural numbers, $E(A) / E(A{\vee} \neg A)$ \eqref{eq:th_defineP}, and can thus be finite or infinite. If $|\mathbb{U}| = |A{\vee} \neg A|$ is infinite, then the evidence for it is also infinite, $E(A{\vee} \neg A) {=} \aleph$ (\ref{sec:th_scale}).  The minimal units of evidence will be $E(u_i) = 1$ and $P(u_i) = 1/\aleph\; (\forall i)$.  The latter is infinitesimal, and $P(u_1) {+} P(u_2) {+} \ldots P(u_\aleph) = 1$ \eqref{eq:th_sum_rule}.  Thus uniformity over an infinite set complies with our sum rule (\ref{sec:th_sum_rule}), but not with the sum rule of the K-axioms and C-axioms (since no real numbers are infinitesimal).   

Hypernatural numbers also solve a related problem that arises in standard measure-theoretic accounts of continuous spaces (\ref{sec:desiderata_specific}).  All finite sets in this context are ``null sets,'' and they are all assigned probability zero \cite{Kapinski2004}.  This is not usually a practical problem, since our interest is typically in comparing infinite sets having finite ratios.  However, we must have $P(u_i) \neq P(u_i {\land} \neg u_i)$, since it is obviously illogical to assign probability zero to both a logical possibility and impossibility.  The problem is solved by allowing for infinitesimal yet variable quantities of evidence. 

If Axiom 4 is used, the formal distinction between discrete and continuous variables vanishes (\ref{sec:ax_scale}).  We always have an infinite set of size $|\mathbb{U}| = \aleph$, even if not all unit propositions (elements) are enumerated.  This overcomes the problems discussed above, but is otherwise inconsequential for probability.  For example, $P(A)$ is the same number whether it is $3/4$ of $1$ or $3/4$ of $\aleph$.  However, the latter is the more apt description in that reality always presents us with infinite possibilities.   

\subsubsection{\label{sec:disc_countably_infinite}The Countably Infinite Continuum}
\vspace{-2mm}
The continuum has long been a source of difficulty in mathematics.  We believe that the solution is to recognize from the start that the infinite is necessary given our lack of knowledge, and that infinities will differ in size.  Once that is accepted, the continuum can be understood as discrete and countable.

Propositions and numbers are inherently discrete and countable.  However, there is no limit to our ability to imagine additional propositions and to invent additional numbers.  This is illustrated by the construction of the surreal numbers, of which the reals and hyperreals are subsets \cite{Knuth1974}.  Starting from 0 and 1, all numbers can be constructed through infinite recursion (a series of discrete steps).  The result is infinities upon infinities, but all countable.  

The standard conclusion that the continuum is uncountable may persist due to the premise that we can always imagine dividing a continuous interval into smaller and smaller intervals, never reaching an interval that is indivisible and therefore discrete.  Our solution is consistent with the premise (of infinite divisibility) but not the conclusion.   The initial interval can only be assigned an arbitrary number by convention (\ref{sec:ax_scale}).  This is often taken to be finite, but we assign a countably infinite number $\aleph$, since the set of intervening possibilities is believed to be infinite.  The interval can then be partitioned into any finite number of sub-intervals, and yet this always leaves us with each sub-interval being infinite, $\aleph/n\; \forall n{\in}\mathbb{N}$, and thus divisible.  This is consistent with our intuitive understanding of continuity as infinite divisibility, yet it preserves countability, since partition into an infinite number of parts, $\aleph$, yields an indivisible unit interval, $\aleph/\aleph$.

The proposed existence of uncountably infinite sets in mathematics, such as $\mathbb{R}$, rests on the Axiom of Choice.  It leads to absurd conclusions and has always been controversial.  We reject its application to infinite sets.  If there actually were some uncountable quantity, we would see it as outside the realm of logic and mathematics.  We also reject the common notion that the continuum of physical space is identical to the continuum of real numbers, as discussed below (\ref{sec:disc_numeric}).  

\subsection{\label{sec:disc_multiple_parameters}Multiple Parameters}
In addition to the problem of infinite sets, uniformity has also been questioned on the grounds that uniformity over one parameter is typically inconsistent with uniformity over another \cite{Berger1985,bernardo1994,Kass1996}. There often appears to be no way of determining which parameter should be uniform, and thus no way to identify $\mathbb{U}$ as objectively correct.

For example, uniformity over angle corresponds to circular symmetry, which is as important and intuitive as any symmetry in nature.  However, nature cannot consist of only one angle, but must feature at least three. A simple ontological model would be `$\mathbf{O} = 3\ particles\ exist$' \eqref{eq:O}, which implies the existence of three angles in a triangle ($\theta_1,\theta_2,\theta_3$). If we had uniformity over one of these angles, we must have non-uniformity over the other two. But that cannot be correct for this model, since it does not quantitatively discriminate between the angles.  

The solution is uniformity over all combinations of the three angles.  However, there are different types of angles in a triangle.  If we have uniformity over all combinations of inner angles (at the vertices) then we cannot have uniformity over all combinations of central angles (at the centroid), and vice versa.  We will show in future work that there is a necessary symmetry that determines the U-space given only the existence of three particles.  

There has not been a general method of finding U-space in probability theory, but the same type of problem was solved in geometry by the Euclidean distance metric. The geometric challenge is eased by the fact that our visual systems readily perceive two dimensions at once.  To better appreciate the problem posed by multiple dimensions, imagine that there is a distance in at least two dimensions, but we have only one-dimensional vision. A person can only perceive and measure the projection of the distance onto a line, and furthermore, each person measures its projection onto a distinct line.  The problem for these observers is that a uniform measure over one line appears inconsistent with a uniform measure over another.

Were our visual abilities so limited, this inconsistency could have caused the field of geometry to become split, as probability theory has (\ref{sec:introduction}), into objectivists and subjectivists, with the subjectivists casting doubt on uniformity as a valid principle, and being skeptical that there is any objectively correct way of measuring distance.  

The specific challenge is to reconcile distance measures $r,x,y$, where $r(a,b)$ is the distance between points $a$ and $b$, and $x(a_x,b_x)$ and $y(a_y,b_y)$ are the distances associated with the orthographic projection of these two points onto two orthogonal axes. The solution is found by noticing that, whereas distance is a measure based on uniformity over a line, angle $\theta$ is a measure based on uniformity over a plane, and it is a function of $x$ and $y$ but not $r$, $\theta = arctan(y/x)$.  Only the Euclidean distance metric, $r = \sqrt{x^2+y^2}$, exhibits rotational invariance and is therefore consistent with uniformity over $\theta$. 

We claim that for any model \textbf{O}, there must be a U-space analogous to a Euclidean space, and a corresponding measure of evidence $E(A{\land}B{\land}C{\land}\ldots)$ analogous to the Euclidean metric \eqref{eq:O}. 

\subsection{\label{sec:disc_absolute_space} Absolute versus Logically Implied Space}
Our axioms imply that the challenge of identifying probabilities is really just the challenge of determining the set of possibilities $\mathbb{U}$ (U-space) given some model of reality.  One of the greatest obstacles is the common assumption that we already know what is possible, and that it corresponds either to ``absolute space,'' or as discussed below, to the numeric domain of a variable (\ref{sec:disc_numeric}).

Given an ontological model \textbf{O} `\textit{there exists a standard deck of $52$ cards},' we can logically derive every possible combination of cards.  There is obviously no other way to identify these possibilities.  Yet the space of possibilities (phase space) in both classical and quantum mechanics is not derived from any such model of what exists.  Rather we simply assume the ``absolute space'' of Newton (which Einstein modified but retained \cite{Norton1993,Maudlin2012}).  This space is proposed to exist independent of, and to be invariant to, whatever matter may be in it, and whatever a local observer may know about that matter.  For example, the possible locations of a particle in Newton's space are the same regardless of whether it is the only particle in the universe, or whether it is 1 of 3 (forming a triangle). 

As famously pointed out by Leibniz and Mach, absolute space has no logical justification, and little or no empirical justification.\footnote{See \cite{Maudlin2012,Barbour2020}. Although Newton's absolute space has been used in models that accurately predict our observations, countless spaces are capable of making the same predictions. For example, the observations explained by Einstein's ``curved spacetime'' can just as well be explained using Newton's ``linear space and time,'' since Einstein's general relativity can be implemented in either \cite{Reichenbach1958}. There can be no empirical evidence favoring one space over the other. Among models that make the same predictions, only non-empirical evidence can favor one over another.} We define an ``absolute space'' to be any space assumed rather than logically deduced from the matter in a model.  In future work we will show how space can instead be deduced from a minimal geometric model of matter.  In contrast to Newton's absolute space, if only 1 particle exists in this model, it can be in only 1 possible location, whereas if 3 particles exist, the space is the infinite set of possible triangles. 

\subsection{\label{sec:disc_numeric}Reification of Numbers}
Even if one discards belief in the existence of absolute space, there is a deeply ingrained view in physics that the space of possibilities is a space of numbers.  Whether intentional or not, scientists often treat numbers in the same way they treat material things they believe to exist.  This reification has profound mathematical consequences that we view as harmful.  

In our view the elements of a set should be entities that exist in a model of reality (at least as possibilities).  Yet the sets at the foundation of mathematics, including measure theory, are sets of purely mathematical objects.  Indeed, Von Neumann provided the widely accepted definition of a number as a set of empty sets (for example, $3 \equiv \{0,1,2\}$, where $0 \equiv \varnothing$, $1\equiv\{\varnothing\}$, and $2 \equiv \{\varnothing,\{\varnothing\}\}$).  In contrast, we define a number not as a set but as only the measure of a set \eqref{eq:three}.  Since we see all sets as countable, including the continuum (\ref{sec:disc_countably_infinite}), we do not distinguish the measure of a set from its cardinality.

Reification of numbers may be more deeply entrenched in physics than the assumption of absolute space.  For example, the physical distance between two bodies is often taken to actually \textit{be} a real number $r{\in}\mathbb{R}^+$.  From this identity it follows that $\mathbb{R}^+$ is the set of possible distances. Although that is permissible with respect to our axioms, since one is free to define model $\mathbf{O}$ as one likes, it is poorly justified with respect to ontology.  Furthermore, uniformity over $\mathbb{R}^+$ conflicts with widespread and well founded beliefs, and it is therefore not accepted as a probability distribution over distance \cite{Berger1985,bernardo1994,Kass1996,Jaynes2003}.  The solution in our view is to revise $\mathbf{O}$ to better reflect actual beliefs.  We argue that a distance is only the measure of a set of actual or possible locations, which taken together constitute a line segment.  Therefore the set of possible distances has a one-to-one mapping to the set of possible line segments.  Determining this set is not as trivial as it may appear (\ref{sec:disc_Uspace}).  

It is relatively well known that there is no uniquely correct number for a single distance considered in isolation, since it cannot be compared to any other distance. Coordinate systems are ultimately arbitrary, and therefore probabilities must be invariant to coordinate transformations \cite{Jeffreys1961,Berger1985,bernardo1994,Jaynes2003}.  But this just means that the space of possibilities does not correspond to a \textit{particular set of numbers}.  We are making the much stronger claim that space does not correspond to \textit{any set of numbers}.  

This may appear radical, but it follows from the fact that we reserve the term `\textit{numbers}' for measures \eqref{eq:three}, whereas the numerals of a coordinate system are only a set of labels for a space.  We all recognize that the numbers associated with the possible outcomes in a game of dice or cards are functioning only as arbitrary labels (and would be better termed `\textit{numerals}').  They are not being used as measures.  But many problems in physics concern geometry, where the space of possibilities is topological and metric. In such cases numbers serve both as measures and labels, and the two can easily become confused. For example, a possibility $u_i{\in}\mathbb{U}$ could be `\textit{a particle lies at $64.2$ degrees from the x-axis}.'  As an angle, this number measures the size of a set of directions. But with respect to $P(u_i {=} 64.2)$, $64.2$ is merely an arbitrary label for a direction, just as \textit{seven of diamonds} is an arbitrary label for a card.  This label may be convenient, but any other could be used, and it need not be numeric (\ref{sec:disc_angle_ratio}).

\subsection{\label{sec:disc_Uspace}Identification of U-space}
Our axioms reduce the problem of deriving probabilities to the problem of deducing possibilities $\mathbb{U}$ from model \textbf{O} \eqref{eq:O}. This can require considerable effort even after overcoming the conceptual obstacles discussed above (\ref{sec:disc_multiple_parameters} - \ref{sec:disc_numeric}).  Here we consider two examples. 
\vspace{-3mm}
\subsubsection{\label{sec:disc_scale_parameter}A Scale Parameter}
\vspace{-2mm}
The ``non-informative prior probability''\footnote{We prefer the phrase ``non-empirical prior.''  This probability necessarily depends on information we have about the numeric representation of a variable, such as positive real numbers versus their logarithms, but it is not dependent on any empirical information we may have about the variable.} of a distance or other scale parameter has been an important and controversial topic.  Since we use $r{\in}\mathbb{R}^+$ as our measure, the state space appears to be $\mathbb{R}^+$ (\ref{sec:disc_numeric}), yet uniformity over $\mathbb{R}^+$ contradicts strongly held beliefs and is widely rejected.  

Harold Jeffreys provided the most widely accepted solution, now known as ``Jeffreys prior over scale.''  He argued that $P(r) = dr/r$, since this is the only probability assignment that is invariant to multiplicative transformation of scale \cite{Jeffreys1932,Jeffreys1961,Jaynes2003}.  But $P(r) = dr/r$ is widely viewed with suspicion on the grounds that it is not consistent with the K-axioms \cite{Berger1985,bernardo1994,Kass1996,Irony1997}.  However, the inconsistency is based on the erroneous assumption that $\mathbb{R}^+$ is the set of possibilities (\ref{sec:disc_numeric}, \ref{sec:disc_angle_ratio}). We intend in future work to derive the possibilities $\mathbb{U}$ and to show how $P(r) = dr/r$ arises from uniformity.  

\subsubsection{\label{sec:disc_angle_ratio}Angle and Ratio}
\vspace{-2mm}
What we often refer to as ``a distance'' or other scale is often a ratio (i.e., it has units), and therefore the U-space consists of at least two dimensions. As a simple but rather contrived example, we can consider ratio $q=y/x$ in a model similar to that introduced above (\ref{sec:disc_multiple_parameters}).  A distance $r{\in}\mathbb{R}^+$ is projected onto arbitrarily chosen orthogonal axes centered on an informative (non-arbitrary) origin.  Therefore $x,y{\in}\mathbb{R}^+$ are distances along these axes, and $r^2 = x^2{+}y^2$.  We would like to identify the set $\mathbb{U}$ of possible ratios $q = y/x$. 

Like the Euclidean metric (\ref{sec:disc_multiple_parameters}), $\mathbb{U}$ can be found from circular symmetry, meaning uniformity over angle $\theta = arctan(y/x)$ (which is the basis of the Gaussian form \cite{Jaynes2003,Fiorillo2021obj}).  However, this does not mean that $\mathbb{U}$ is the set of real numbers that constitutes the domain of $\theta$.  Rather it means that we will have uniformity if we partition the plane into $\mathfrak{n}$ tranches of equal angle as measured from the origin.  If we choose $\mathfrak{n}{=}90$ tranches, each will span $1$ degree.  But since angle is infinitely divisible, to find all possibilities $u_1,u_2\ldots$ we need to partition this space into an infinite number of tranches $\mathfrak{n{\in}\aleph}$, so that each $u_i$ spans $90/\aleph$ degrees (\ref{sec:ax_scale}, \ref{sec:disc_infinite_solution}). The result is $P(u_i) = 1/\aleph\ (\forall i)$.  

We naturally say that we have ``uniformity over angle $\theta$,'' but it should be recognized that $\theta$ is only an arbitrary label for a possibility $u_i$, not the possibility itself (\ref{sec:disc_numeric}). We can have any number of labels for a single possibility, another being $q = y/x = tan(\theta)$.  The labels make no difference to the true set of possibilities $\mathbb{U}$, over which we always have uniformity.  However, we lose uniformity if we group these possibilities into subsets based on equal intervals of $q$ (\ref{sec:spaces}, \ref{sec:disc_multiple_parameters}). 

We demonstrated previously that $P(q{<}1) = P(q{>}1)$, so that $q{=}1$ is the median of the distribution, not just in this case but for a large set of statistical models that includes all those famous enough to have names \cite{Fiorillo2021obj}. The present axioms indicate that it is not only probability, but space itself, that is symmetrically distributed around $q{=}1$. Thus there are the same number of possibilities for $0{<}q{\leq} 1$ as for $1{\leq}q{<} \infty$.  The fact that most numbers in $\mathbb{R}^+$ are greater than $1$ is irrelevant, because these are not alternative possibilities for $q$, but alternative measures for $q$ (\ref{sec:disc_numeric}).

\section{\label{sec:conclusion}Conclusion}
Probability theory has been a remarkably contentious topic for at least the last hundred years, and we do not expect that to change soon.  We nonetheless hope the present work might help to persuade frequentists and subjectivist Bayesians that probabilities can be objectively correct measures of evidence in cases where one is able to precisely define all relevant propositions.  Our demonstration that uniformity must underlie all rigorously justified probabilities should help to overcome technical if not philosophical disputes.

We began by noting that there has not been an established measure of evidence.  We can rephrase this by stating that there has not been an established and generally applicable method of determining what is possible. Our axioms assert that if we can specify exactly what is possible, given some model, then the evidence for every possibility must be equal. Furthermore, since we do not distinguish a `set of possibilities' from a `space of possibilities,' or simply a `space,' the measure of evidence is only the measure of space.

\section*{Acknowledgments}
{The authors declare they have no financial or other conflicts of interest.}

\bibliography{Axioms}
\end{document}